\documentclass[smallextended]{svjour3}       % onecolumn (second format)
\smartqed  % flush right qed marks, e.g. at end of proof
\usepackage{graphicx}
\journalname{Numer.\ Math.}
\def\Im{\hbox{\rm Im\kern .7pt}}
\def\Re{\hbox{\rm Re\kern .7pt}}
\def\arg{\hbox{\rm arg\kern .2pt}}
\def\twopi{2\kern.2 pt\pi}

\def\zc{z_c^{}}
\def\Rn{R_n}
\def\MATLAB{\hbox{MATLAB}}
\def\inz{[\kern .5pt 0,1]}
\def\inzz{[\kern .5pt 0.01,1]}
\def\inzm{[-1,0\kern .5pt]}
\def\sx{{\sqrt x\kern 1pt}}
\def\ve{\varepsilon}
\def\uE{u_E^{}}
\def\uG{u_\Gamma^{}}
\def\Imin{I_{\min{}}}
\def\sech{\hbox{sech}}
\def\pih{{\textstyle{\pi\over 2}}}
\def\half{{\textstyle{1\over 2}}}

\titlerunning{Exponential node clustering at singularities}
\title{Exponential node clustering at singularities for
rational approximation, quadrature, and PDE\kern .5pt s}

\author{Lloyd N. Trefethen \and Yuji Nakatsukasa \and J. A. C. Weideman}

\institute{L. N. Trefethen \at
Mathematical Institute\\
University of Oxford\\
Oxford OX2 6GG, UK\\
\email{trefethen@maths.ox.ac.uk} \\
\and
Y. Nakatsukasa\at
Mathematical Institute\\
University of Oxford\\
Oxford OX2 6GG, UK\\
\email{nakatsukasa@maths.ox.ac.uk} \\
\and
J. A. C. Weideman\at
Applied Mathematics\\
Stellenbosch University\\
Private Bag X1\\
Matieland 7602, South Africa\\
\email{weideman@sun.ac.za} \\
}

\date{Received: date / Accepted: date}

\usepackage{graphics,graphicx}
\usepackage{upquote,verbatim,amssymb}
\def\C{{\bf C}}

\begin{document}

\maketitle

\begin{abstract}
Rational approximations of functions with singularities can
converge at a root-exponential rate if the poles are exponentially
clustered.  We begin by reviewing this effect in minimax,
least-squares, and AAA approximations on intervals and complex
domains, conformal mapping, and the numerical solution of Laplace,
Helmholtz, and biharmonic equations by the ``lightning'' method.
Extensive and wide-ranging numerical experiments are involved.
We then present further experiments showing that in all of these
applications, it is advantageous to use exponential clustering
whose density on a logarithmic scale is not uniform but tapers
off linearly to zero near the singularity.  We give a theoretical
explanation of the tapering effect based on the Hermite contour
integral and potential theory, showing that 
tapering doubles the rate of convergence.  Finally we show that related
mathematics applies to the relationship
between exponential (not tapered) and doubly exponential (tapered)
quadrature formulas.  Here it is the Gauss--Takahasi--Mori contour
integral that comes into play.
\end{abstract}

\section{\label{sec-intro}Introduction}
Analytic functions can be approximated by polynomials with
exponential convergence, i.e., $\|f-p_n\| = O( \exp(-Cn))$ for
some $C>0$ as $n\to\infty$.  Here $n$ is the polynomial degree
and $\|\cdot\|$ is the $\infty$-norm on an approximation domain
$E$, which may be a closed interval of the real axis or more
generally a simply connected compact set in the complex plane.
This result is due to Runge~\cite{runge,walsh} and explains the
exponential convergence of many numerical methods when applied
to analytic functions, including Gauss and Clenshaw--Curtis
quadrature~\cite{gaussCC,atap} and spectral methods for ordinary
and partial differential equations~\cite{smim,series}.  It is
also the mathematical basis of Chebfun~\cite{chebfun}.

If $f$ is not analytic in a neighborhood of $E$, then
Bernstein showed in 1912 that exponential convergence of
polynomial approximations is impossible~\cite{bernstein,atap}.
Bernstein also showed that in approximation of functions with derivative
discontinuities such as $f(x) = |x|$ on $[-1,1]$, polynomials can converge
no faster than $O(n^{-1})$~\cite{b14}.  Now from the beginning,
going back to Chebyshev in the mid-19th century, approximation
theorists had investigated approximation by rational functions as
well as polynomials.  Yet it was not until fifty years after these
works by Bernstein that it was realized that for this problem
of approximating $|x|$ on $[-1,1]$, rational functions
can achieve the much faster rate of {\em root-exponential
convergence,} that is, $\|f-r_n\| = O( \exp(-C\sqrt n\kern
1pt))$ for some $C>0$.  This result was published by
Newman in 1964~\cite{newman}, who also showed that faster
convergence is not possible.  With hindsight,
it can be seen that the root-exponential effect was implicit in
the results of Chebyshev's student Zolotarev nearly a century
earlier~\cite{gonZ,nf,stahl93,zol}, but this was not noticed.

Newman's theorem has been a great stimulus
to further research in rational approximation
theory~\cite{gonZ,gon67,gonchar,levsaff,safftotik,stahlgeneral,stahl93,vyach}.
It has not, however, had much impact on scientific computing
until very recently with the discovery that it can be the basis of
root-exponentially converging numerical methods for the solution
of partial differential equations (PDE\kern .3pt s) in domains
with corner singularities~\cite{stokes,conf,lightning,PNAS,laplace}.
The aim of this paper is to contribute to building the bridge
between approximation theory and numerical computation.

In particular, we shall focus on the key feature that gives
rational approximations their power: the exponential clustering of
poles near singularities.  (The zeros are also
exponentially clustered, typically interlacing the poles, with the
alternating pole-zero configuration serving as proxy for a branch cut.)
This has been a feature of the theory since Newman's explicit
construction.  Our aim is, first, to show how widespread this
effect is, not only with minimax approximations (i.e., optimal in
the $\infty$-norm), the focus of most theoretical studies, but
also for other kinds of approximations that may be more useful
in computation.  Section~\ref{sec2} explores this effect in a
wide range of applications.

In section~\ref{sec3} we turn to a new contribution of this
paper, the observation that good approximations tend to
make use of poles which, although exponentially clustered,
have a density on a logarithmic scale that tapers to zero at
the endpoint.  Specifically, the distances of the clustered
poles to the singularity appear equally spaced when the log of
the distance is plotted against the square root of the index.
We show experimentally that this scaling appears not just with
minimax approximations but more generally.

To explain this effect, we begin with a review in section~\ref{sec4}
of the Hermite contour integral, which is the basis of the application of
potential theory in approximation.  We show how
this leads to the idea of condenser capacity for
the analysis of rational approximation of analytic functions.
Section~\ref{sec5} then turns to functions with singularities
and explains the tapering effect.
In this case the condenser is short-circuited, and it is
not possible to estimate the Hermite integral by considering
the $\infty$-norm of the factors of its integrand, but the
$1$-norm gives the required results.  Analysis of a model problem
shows how the tapered exponential clustering of
poles enables better overall resolution, potentially doubling
the rate of convergence.  These arguments are
related to those developed in the
theoretical approximation theory literature by Stahl and
others~\cite{stahlgeneral,stahl93,stahl}, but we believe that
section~\ref{sec5} of this paper is the first to connect this
theory with numerical analysis.

Finally in section~\ref{sec6} we turn to a different
problem, the quadrature of functions with endpoint
singularities on $[-1,1]$.  Here the famous methods are
the exponential (tanh) and double exponential (tanh-sinh)
formulas~\cite{haber,IMT,mori,ms01,oms,sugihara,tm73,tm74,tanaka}.
Making use of the link to another contour integral formula, the
Gauss--Takahasi--Mori integral~\cite{gauss,gaussCC,tm71}, we show
that the distinction between straight and tapered exponential
clustering arises here too.

Throughout the paper, $\Rn$ denotes the set of rational functions
of {\em degree~$n$}, that is, functions that can be written as
$r(x) = p(x)/q(x)$ where $p$ and $q$ are polynomials of degree~$n$.
The norm $\|\cdot\|$ is the $\infty$-norm on $E\kern 1pt$,
but, as mentioned above, other measures will come into play in sections~\ref{sec5}
and~\ref{sec6}, and indeed, a theme of our discussion is that
certain aspects of rational approximation are often concealed by
too much focus on the $\infty$-norm.

The numerical experiments in this paper are a major part
of the contribution; we are not aware of comparably detailed studies
elsewhere in the literature.
Our emphasis is on the results, not the algorithms, but our
numerical methods are briefly summarized in the discussion section at the end.

\section{\label{sec2}Root-exponential convergence and exponential clustering of poles}
\begin{figure}
\vskip 1.4em
\begin{center}
\includegraphics[scale=.99]{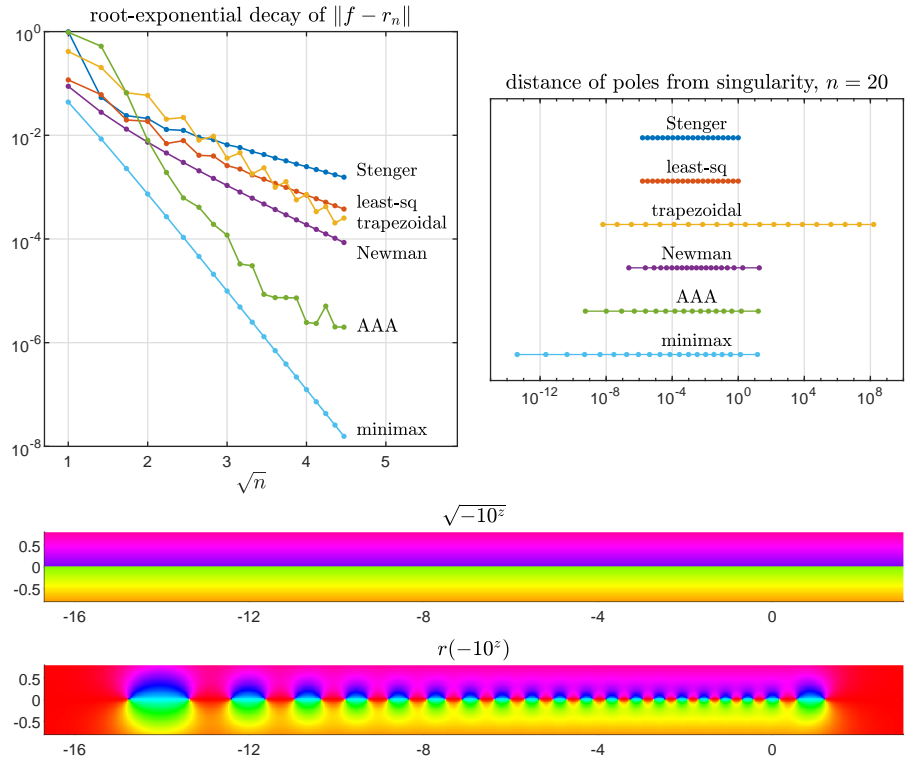}
\end{center}
\caption{\label{fig1}Root-exponential convergence of six kinds
of degree $n$ rational approximations
of $f(x) = \sqrt x$ on $\inz$ as $n\to\infty$.  
On the upper-left, the asymptotically straight
lines on this log scale with $\sqrt n$ on the horizontal
axis (except for AAA) show the root-exponential effect.
On the upper-right,
the distances of the poles in $(-\infty,0)$ from the singularity 
at $x=0$ show the exponential clustering.  Below, phase
portraits in the complex plane
of the square root function (the standard branch)
and its degree 20 minimax approximation on $\inz$, after an exponential
change of variables, show how a branch cut is approximated
by interlacing exponentially clustered poles and zeros. 
Red before yellow going counterclockwise indicates a zero,
and yellow before red indicates a pole.
We use $10^z$ instead of $e^z$ to enable comparison with the axis labels
in the images above.}
\end{figure}

\begin{figure}
\begin{center}
\vskip 1em
\includegraphics[scale=.99]{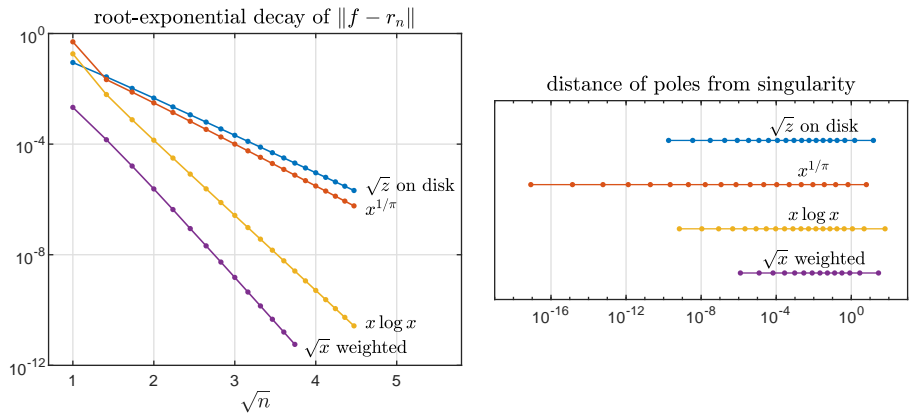}
\end{center}
\caption{\label{fig2}Four more minimax approximations, showing
the same root-exponential convergence and exponential
clustering of poles as in Fig.~{\rm\ref{fig1}}.  Two involve
the functions $x^{1/\pi}$ and $x\log x$ on $\inz$, one involves
$x$ on $\inz$ but with the $\infty$-norm weighted by $x$, and one
involves $\sqrt z$ on the disk about $\half$ of radius $\half$.
In the right image, $n$ takes its final value from the left
image for each problem, $14$ for the weighted
approximation and $20$ for the other cases.}
\end{figure}

In this section we explore the convergence
of a variety of rational approximations to analytic functions
with boundary branch point singularities.  Our starting point
is Fig.~\ref{fig1}, which presents results for six kinds
of approximations of $f(x) =\sqrt x$ on $\inz$ by rational
functions of degrees $1\le n \le 20$.  (By the substitution $x =
t^2$, this is equivalent to Newman's problem of approximation
of $|t|$ on $[-1,1]$.)  The choice of $f$ is not special; as
we shall illustrate in Figs.~\ref{fig2} and~\ref{fig4},
other functions with endpoint singularities give similar results.

First, the big picture.  The upper-left image of the figure shows
$\infty$-norm errors $\|f-r_n\|$ plotted on a log scale as
functions not of $n$ but of $\sqrt n$.  With the exception of the
erratic case labeled AAA, all the curves plainly approach straight lines
as $n\to\infty$: root-exponential convergence.  (The shapes would
be parabolas if we plotted against $n$.)  The upper-right image shows
the absolute values of the 20 poles for the approximations with
$n=20$, that is, their distances from the singularity at $x=0$.
On this logarithmic scale the poles are smoothly distributed:
exponential clustering.  This clustering is further shown
in the lower images, for the approximation labeled minimax,
by a phase portrait~\cite{wegert} of the square root function
(the standard branch) and its degree 20 rational approximation
after an exponential change of variables.

The top four approximations have preassigned poles, making the
approximation problems linear; indeed the Stenger, trapezoidal,
and Newman approximations are given by explicit formulas.
The AAA and minimax approximations are nonlinear,
with poles determined during the computation.  
Although it is tempting to rank these candidates from
worst at the top to best at the bottom (the minimax approximation
is best by definition), this is not the point.
All these approximations converge root-exponentially, and the
differences in efficiency among them amount to constant factors of
order 10, which can in fact be improved in most cases by
introducing a scaling parameter or two.  In particular,
minimax and other nonlinear approximations can approximately
double the rate of convergence of
the linear approximations~\cite{rakh}.
All these approximations can achieve accuracy $10^{-6}$
with degrees $n\approx 100$, whereas with polynomials one
needs $n=\hbox{140,085}$.

We comment now on the individual approximations of
Fig.~\ref{fig1}.  The Newman approximation comes
from the explicit formula presented in his four-page
paper~\cite{newman}.  The approximation is\/ $r(x) = \sx (\kern
.7pt p(\sx)-p(-\sx))/(\kern .7pt p(\sx)+p(-\sx))$, where $p(t)
= \prod_{k=0}^{2n-1} (t+\xi^k)$ and $\xi = \exp(-1/\sqrt{\kern
.5pt 2n}\kern 1pt )\kern .7pt$; this can be shown to be a rational
function in $x$ of degree $n$.  The asymptotic convergence rate is
$\exp(-\sqrt{\kern .5pt 2\kern .3pt n}\kern 1pt )$~\cite{xz}.  This
can be improved to approximately $\exp(-(\pi/2)\sqrt{\kern .5pt 2\kern
.3pt n}\kern 1pt )$ by defining $\xi = \exp(-(\pi/2)/\sqrt{\kern
.5pt 2\kern .3pt n}\kern 1pt )$, an example of the scaling
parameters mentioned in the last paragraph (these values are
conjectured to be optimal based on numerical experiments).

The trapezoidal approximation
originates with Stenger's investigations
of sinc functions and associated
approximations~\cite{stengersurvey,stenger,stengerbook}.
Following p.~211 of~\cite{atap}, we approximate
$\sqrt x$ by starting from the identity
$\sqrt x = {2\kern .3pt x/\pi} \int_0^\infty (t^2 + x)^{-1}dt$,
which with the change of variables $t = e^s$ becomes
\begin{equation}
\sqrt x = {2\kern .3pt x\over \pi} \int_{-\infty}^\infty {e^s ds\over e^{2s} + x}.
\end{equation}
For $n\ge 1$, we approximate this integral by an equispaced
$n$-point trapezoidal rule with step size $h>0$,
\begin{equation}
r(x) = {2h x\over \pi} \sum_{k =-(n-1)/2}^{(n-1)/2} 
{e^{kh}\over e^{2kh}+x}.
\label{traprule}
\end{equation}
(If $n$ is even, the values of $k$ are half-integers.)  There are
$n$ terms in the sum, so $r$ is a rational function of degree
$n$ with simple poles at the points $p_k = -\exp(2kh)$.
Two sources of error make $r(x)$ differ from $\sqrt x$.
The termination of the sum at $n<\infty$ introduces an
error of the order of $\exp(-nh/2)$, and the finite step
size introduces an error on the order of $\exp(-\pi^2/h)$,
since the integrand is analytic in the strip around the
real $s$-axis of half-width $\pi/2$~\cite[Thm.~5.1]{trap}.
Balancing these errors gives the optimal step size $h \approx
\pi\sqrt{\kern .5pt 2/n}$ and approximation error $\|r - \sx \|
\approx \exp(-\pi\sqrt{n/2}\kern .5pt)$.  Note that the poles
for this approximation cluster at $\infty$ as well as at $0$,
and indeed, it converges root-exponentially not just on $\inz$
but on any interval $[\kern .5pt 0,L]$ with $L>0$.

The derivation by the trapezoidal rule just given explains in
a general way why root-exponential convergence is achievable
for a wide range of problems with endpoint singularities.
With any exponentially graded discretization, there will be
errors associated with finite grid sizes and errors associated
with truncation of an infinite series.  If both sources of error
follow an exponential dependence, then an optimal balance with
step sizes scaling with $1/\sqrt n$ can be expected to lead to a
root-exponential result.  Such effects are familiar in the analysis
of $hp$ discretizations of partial differential equations when
the step sizes $h$ and orders $p$ of multiscale discretizations
are balanced to achieve optimal rates of convergence near
corners~\cite{schwab}.

A drawback of the trapezoidal approximation is that its derivation
depends on the precise spacing of the poles, since it relies on
the property that the trapezoidal rule is exponentially accurate
in this special case~\cite{trap}.  The curves labeled Stenger in
Fig.~\ref{fig1} come from a more flexible alternative approach,
also proposed by Stenger~\cite{stenger}, where we fix $n$
distinct poles $p_k\in (-\infty,0)$, $1\le k \le n$ and $n+1$
interpolation points $x_k\in \inz$, $0\le k \le n$, and then
take~$r$ to be the unique rational function of degree $n$ with
these poles that interpolates $f(x)$ in these points.  The theory
of rational interpolation with preassigned poles was developed by
Walsh~\cite{walsh} and will be discussed in section~\ref{sec4}.
For our problem of approximation on $\inz$ with a singularity at
$x=0$, a good choice is to take $x_0 = 0$ and $x_k = -p_k$ for
$k\ge 1$.  In particular, our {\em Stenger approximant\/}\footnote{Stenger
considered rational approximations of this kind, though not in
this precise setting of a finite interval with just one endpoint
singularity.} is the rational function $r$ resulting from the
choices \begin{equation} -p_k = x_k = \exp(-(k-1) h) , \quad 1\le
k \le n, \label{stenginterp} \end{equation} with $h = O(1/\sqrt
n\kern 1pt)$.  Figure~\ref{fig1} takes $h = \pi /\sqrt n$.

Interpolation is important for theoretical analysis, but for
practical computation, least-squares fitting is more robust
and more accurate, since it does not require knowledge of good
interpolation points.  The least-squares data of Fig.~\ref{fig1}
come from fixing the same exponentially clustered poles as in
(\ref{stenginterp}), but now choosing approximation coefficients
by minimizing the least-squares error $f-r$ on a discretization
of $\inz$ by standard methods (\MATLAB\ backslash).  As always
when discretizing near singularities, we use an exponentially
graded mesh ({\tt logspace(-12,0,2000})), and a weight function
$w(x) = \sqrt x$ is introduced in the discrete least-squares
problem so that it approximates a uniformly weighted problem on
the continuum.  The error curve $r(x) - \sqrt x$ for $x\in \inz$
for this approximation (not shown) approximately equioscillates
between $n+2$ extrema, indicating that it is a reasonable
approximation to the best $L^\infty$ approximation with these
fixed poles.

The minimax data in Fig.~\ref{fig1} correspond to the true
optimal (real) approximations, rational approximations with free poles.
Here the error curve equioscillates between $2\kern .3pt n+2$
extrema~\cite{atap}, and the error is approximately squared;
the asymptotic convergence rate is $\exp(-\pi\sqrt{\kern .5pt
2\kern .3pt n}\kern 1pt )$~\cite{stahl,vyach}.

Computing minimax approximations, however, can be challenging~\cite{minimax},
and on a complex domain they need not even be unique~\cite{gt}.
This brings us to the data in the figure for AAA (adaptive Antoulas--Anderson)
approximation, a fast method of near-best rational approximation
introduced in~\cite{aaa}.  AAA approximation is at
its least robust on real intervals, as reflected in the erratic
data of the figure, but for more complicated problems and in the
complex plane, it is often the most practical method for
rational approximation.

This concludes our discussion of Fig.~\ref{fig1}.  The next
figure, Fig.~\ref{fig2}, illustrates that these effects are not
confined to approximation on a real interval or to the function
$\sqrt x$.  The figure presents data for four further examples of
minimax approximations.  One set of curves shows approximation of
$x^{1/\pi}$ on $\inz$, with the value $1/\pi$ chosen to dispel any
thought that rational exponents might be special.  This problem
requires poles particularly close to the singularity since the
exponent is so small.  Another shows approximation of $x\log x$
on $\inz$.  With a much weaker singularity, this problem shows
higher approximation accuracy.  A third shows approximation of
$\sqrt x$ again, but now it is weighted minimax approximation,
with a weight function $x$ (and the error measured is now the
weighted error, notably smaller than before).  Finally the
fourth set of data shows minimax approximation of $\sqrt z$
on the complex disk $\{z {\kern 1pt:\kern 3pt} |z-\half|< \half\}$.

\begin{figure}
\begin{center}
\vskip 1em
\includegraphics[scale=.95]{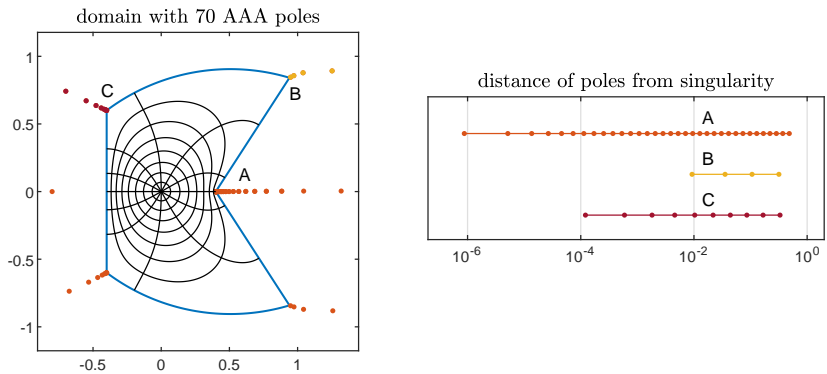}
\end{center}
\caption{\label{fig3}The conformal map of a circular pentagon onto
the unit disk has been computed and then approximated numerically
by a rational function of degree $70$~{\rm\cite{conf,conformal}} by the
AAA algorithm.  The poles cluster exponentially at the corners, where the
map is singular.}
\end{figure}

Figure~\ref{fig3} turns to our first problem of scientific
computing.  Following methods presented in~\cite{conf}
and~\cite{conformal}, a region $E$ in the complex plane bounded
by three line segments and two circular arcs has been conformally
mapped onto the unit disk, and the map has then been approximated
to about eight digits of accuracy by AAA approximation, which
finds a rational function with $n=70$.  This process is entirely
adaptive, based on no a priori information about corners or
singularities, yet it clusters the poles near the corners just
as in Figs.~\ref{fig1} and~\ref{fig2}.\ \ Many poles cluster
at the strong singularity {\sf A} and only a few at the weak
singularity {\sf B}.\ \ Note that the poles lie asymptotically on
the bisectors of the external angles.  This effect is well known
especially from the theory of Pad\'e approximation as worked out
initially by Stahl~\cite{stahl12,suetin}.  Optimal approximations line up their
poles along curves which balance the normal derivatives of a
potential gradient on either side, and evidently the AAA method
comes close enough to optimal for the same effect to appear.

We finish this section with a look at lightning solvers for
PDE\kern .3pt s in two-dimensional domains, introduced in 2019
and applied to date to Laplace~\cite{lightning,PNAS,laplace},
Helmholtz~\cite{PNAS}, and biharmonic equations
(Stokes flow)~\cite{stokes}.  In the basic case of a Laplace problem
$\Delta u = 0$, the idea is to represent the solution on a domain
$E$ as $u(z) \approx \Re r(z)$, the real part of a rational
function with no poles in $E$ that approximates the boundary data
to an accuracy typically of 6--10 digits.  The rational
functions have preassigned poles that cluster exponentially at
the corners, where the solution will normally have singularities~\cite{lehman,wasow},
and the name ``lightning'' alludes to this exploitation of the same
mathematics that makes lightning strike objects at sharp corners.
Coefficients for the solution are found by least-squares fitting,
making this an approximation process of the same structure as
in the least-squares example of Fig.~\ref{fig1}.\ \ The difference
is that the approximations are now applied to give values of
$u(z)$ in the interior of the domain $E$, where it is not known
a priori.  See Fig.~\ref{fig4} for an example on a ``snowflake''
with boundary data $\log|z|$.

\begin{figure}
\begin{center}
\vskip 2em
\includegraphics[scale=.85]{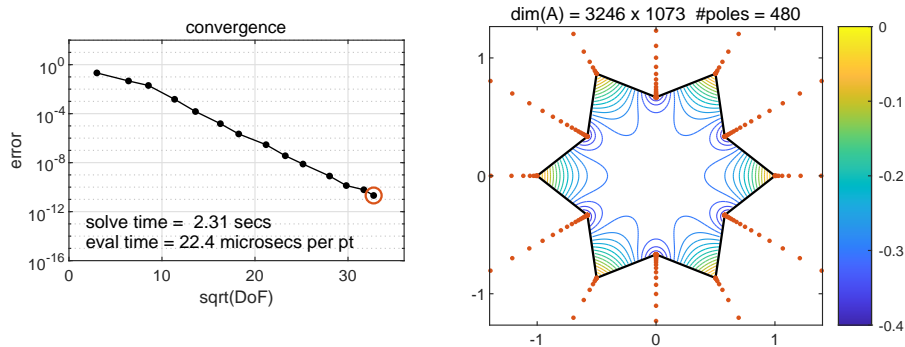}
\end{center}
\caption{\label{fig4}Example of the lightning Laplace 
solver~{\rm\cite{lightning,PNAS}} as implemented
in the code {\tt laplace.m}~{\rm\cite{laplace}}.
For each number of degrees of freedom (DoF),
poles are clustered exponentially near the
$12$ corners of the domain $E$, and the numbers are
increased until a solution to $10$-digit accuracy
is obtained in the form of a rational function
with $480$ poles.  This takes $2.3$ s on a laptop, and subsequent
evaluations take $22$ $\mu$s per point, with the accuracy of each
evaluation guaranteed by the maximum principle.}
\end{figure}

\begin{figure}
\begin{center}
\vskip 10pt
\raisebox{12pt}{\includegraphics[scale=.66]{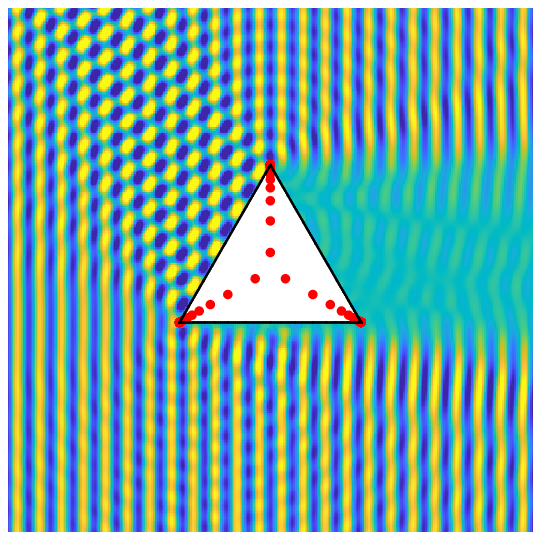}~~~~~~~~~~~~~~~~~}%
\includegraphics[scale=.80]{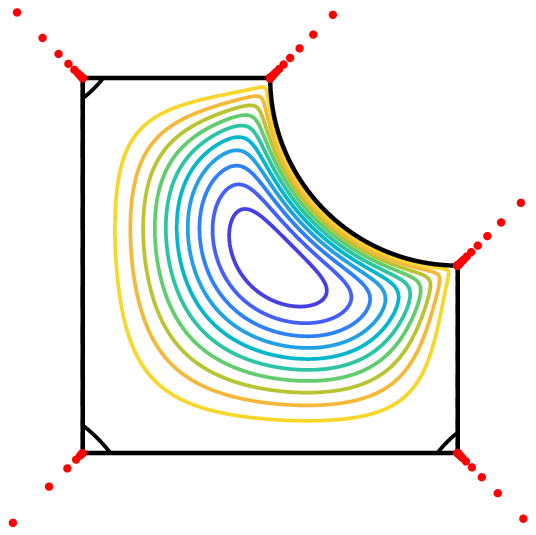}
\vspace{-5pt}
\end{center}
\caption{\label{fig5}Lightning solvers have been generalized to the
two-dimensional Helmholtz (left)~{\rm\cite{PNAS}} and biharmonic
equations (right)~{\rm\cite{stokes}}.
The Helmholtz image shows a plane wave incident from
the left scattered from
a sound-soft equilateral triangle.  The biharmonic image shows
contours of the stream function for Stokes flow in a cavity driven
by a quarter-circular boundary segment rotating at
speed 1 and with zero velocity on the remainder of the boundary.
The black contours in the corners, representing the stream
function value $\psi = 0$, delimit counter-rotating Moffatt vortices.
Tapered exponentially clustered singularities are
used in both computations.}
\end{figure}

Lightning solvers have been generalized to the Helmholtz
equation $\Delta u + k^2u = 0$~\cite{PNAS} and the biharmonic
equation $\Delta^2 u = 0$~\cite{stokes}, as illustrated in
Fig.~\ref{fig5}.  In the Helmholtz case, poles $(z-z_k)^{-1}$
of rational functions become singularities of complex Hankel functions
$H_1(k|z-z_k|)\exp(\pm i \kern .7pt \arg(z-z_k))$, and the biharmonic case is handled by the
Goursat reduction $u(z) = \Im (\kern .5pt  \overline z f(z) +
g(z))$ to a coupled pair of analytic functions $f$ and $g$,
each of which is approximated by its own rational function.
The mathematics of lightning methods for Helmholtz and biharmonic
problems has not yet been worked out fully, and the analysis
given in section~\ref{sec5} applies just to the Laplace case.

Although it is not the purpose of this article to give details
about lightning PDE solvers, they are at the heart of our motivation.
Usually in approximation theory, minimax approximations are
investigated as an end in themselves, and the locations of
their poles may be examined as an outgrowth of this process;
a magnificent example is~\cite{stahl94}.  Here, the order is
reversed.  Our aim is to exploit an understanding of how poles
cluster to construct approximations on the fly to solve
problems of scientific computing.

\section{\label{sec3}Tapered exponential clustering}
In the last section, 13 plots were presented of the distances
of poles to singularities on a log scale, the right-hand images
of Figs.~\ref{fig1}, \ref{fig2}, and~\ref{fig3}.  All showed
exponential clustering, and all but three showed a further effect
which we call {\em tapered exponential clustering\/\kern .4pt},
the main subject of the rest of this paper: on the log scale, the
spacing of the poles grows sparser near the singularity.  This was
also colorfully evident in the phase portrait at the bottom
of Fig.~\ref{fig1}.  The three exceptions were the Stenger,
least-squares, and trapezoidal approximations of Fig.~\ref{fig1},
all of which are based on poles preassigned with strictly uniform
exponential clustering.  These examples illustrate that tapering
of the pole distribution is not necessary for root-exponential
convergence.
A fourth set of data in Fig.~\ref{fig1} also involves preassigned
poles, the Newman data, and some tapering is apparent in this case.

\begin{figure}
\begin{center}
\vskip 15pt
\includegraphics[scale=1]{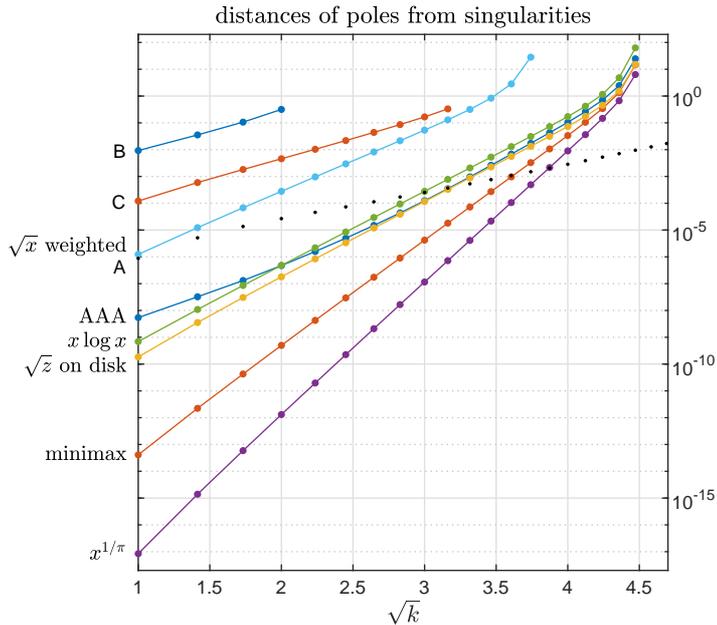}
\vspace{-8pt}
\end{center}
\caption{\label{fig6}Tapered exponential clustering of poles
near singularities for the nine examples with free poles
from Figs.~{\rm\ref{fig1}--\ref{fig3}} of the last section.
The crucial feature is that the curves appear straight with
this horizontal
axis marking $\sqrt k$ rather than $k$, where $\{d_k\}$ are
the sorted distances of the poles from the singularities.
The data for the poles at vertex {\sf A} of 
Fig.~{\rm \ref{fig3}}
have been deemphasized to diminish clutter (black dots),
since they lie at such a different slope from the others.}
\end{figure}

Figure~\ref{fig6} shows the nine remaining examples of exponential
clustering of poles from Figs.~\ref{fig1}--\ref{fig3},
the ones with free poles,
presenting the distances $\{d_k\}$ of the
poles from their nearest singularities on a log scale.  What is
immediately apparent is that all the curves look straight
for smaller values of $k$.  Note that five of them stop at $n=20$,
one at $n=14$, and the remaining three, from the approximation
of a conformal map of Fig.~\ref{fig3}, at different values
determined adaptively by the AAA algorithm.

Yet the horizontal axis in Fig.~\ref{fig6} is not $k$ but $\sqrt k$.
Plotted against $k$ (not shown), the data would look completely
different.  Evidently in a wide range of rational
approximations, both
best and near-best, the distances $\{d_k\}$ of poles to singularities
is well approximated by the formula
\begin{equation}
\log d_k \approx \alpha + \sigma \sqrt k
\label{model1}
\end{equation}
for some constants $\alpha$ and $\sigma$, that is,
\begin{equation}
d_k \approx \beta\exp(\sigma \sqrt k \kern 1pt)
\label{model2}
\end{equation}
for some $\beta$ and $\sigma$.

\begin{figure}
\vskip 15pt
\begin{center}
\includegraphics[scale=1.2]{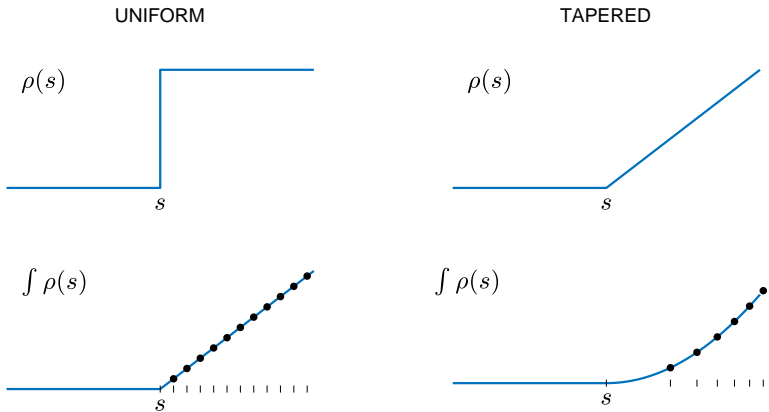}
\vspace{-3pt}
\end{center}
\caption{\label{relu}The algebra of exponential clustering.
With respect to the variable $s = \log d$, where $d$ is the distance to the
singularity, the simplest exponential clustering of poles
would have uniform density $\rho(s)$ 
down to a certain value and then cut off abruptly
(left column).  A tapered distribution cuts off linearly instead (right column),
resulting in poles exponentially clustered in
the $\sqrt k$ fashion seen in Fig.~{\rm \ref{fig6}}.}
\end{figure}

To make sense of the $\sqrt k$ scaling, let us remove the exponential
from the problem by defining a distance variable $s= \log d$,
thereby transplanting an interval such as $d\in \inz$ to $s\in
(-\infty, 0\kern .5pt]$.  We ask, what can be said of the density
$\rho(s)$ of poles with respect to $s\kern .7pt$?  If $\rho(s)$
were constant, this would correspond to a uniform exponential
distribution of poles, requiring an infinite number of poles
since $s$ goes to $-\infty$.  So some kind of cutoff of $\rho(s)$
to $0$ must occur as $s\to-\infty$.  An abrupt cutoff, as with
the Stenger, trapezoidal, and least-squares distributions of
Fig.~\ref{fig1}, leads to a linear cumulative distribution,
as shown in the left column of Fig.~\ref{relu}.  By contrast, a linear
cutoff gives a quadratic cumulative distribution, as shown in the
right column, and when this
is inverted, the result is the $\sqrt k$ distribution we have
observed.

Thus the straight lines of Fig.~\ref{fig6} can be explained
if pole density functions $\rho(s)$ for good rational 
approximations tend to take the form 
sketched in the upper-right of
Fig.~\ref{relu}.  (Aficionados of deep learning may call this
the ``ReLU'' shape.)  In section~\ref{sec5} 
we will explain why this is the case and continue the story
of Fig.~\ref{relu} in Fig.~\ref{relupot}.

We have not presented data in this section for lightning PDE solutions,
but it was in this context
that we first became aware of the importance of tapered
exponential clustering.
In the course of the work leading to~\cite{lightning},
the first author noticed that
although straight exponential spacing of preassigned poles gave
root-exponential convergence, better efficiency could be achieved
if the resulting approximations
were re-approximated a second time by the AAA algorithm.
On examination it was found that the AAA
approximations had poles in a
tapered distribution, just like cases {\sf A}--{\sf C} 
of Fig.~\ref{fig6}.  The model
(\ref{model1})--(\ref{model2}) was developed empirically in this
context, with $\sigma\approx 4$ found to be an effective choice.
This became
the formula for preassignment of poles in the lightning Laplace
software~\cite{laplace}, where it improved the overall speed
by a good factor, and it appears as equation (3.6) in~\cite{lightning}.

\section{\label{sec4}Hermite integral formula and potential theory}
The basic tool for estimating accuracy of rational
approximations is the Hermite integral formula~\cite{levsaff,walsh}.
In this section we review how this
formula leads to the use of potential theory~\cite{ransford}, and
in particular the quantity known as the
condenser capacity, for approximations of analytic functions.
Building on the work of Walsh~\cite{walsh}, these ideas began to
be developed by Gonchar and Rakhmanov in the Soviet Union not long
after the appearance of Newman's paper~\cite{gon67,gonchar}.

The following statement is adapted from Thm.~8.2 of~\cite{walsh}.

\begin{theorem}
Let\/ $\Omega$ be a simply connected domain in $\C$ bounded
by a closed curve $\Gamma$, and let $f$ be analytic in 
$\Omega$ and extend continuously to the boundary.
Let distinct interpolation points $x_0,\dots, x_n\in \Omega$
and poles $p_1,\dots,p_n$ anywhere in
the complex plane be given.
Let $r$ be the unique degree $n$ rational function
with simple poles at $\{p_k\}$ that
interpolates $f$ at $\{x_k\}$.  Then for any $x\in\Omega$,
\begin{equation}
f(x) - r(x) = {1\over 2\pi i} \int_\Gamma {\phi(x)\over \phi(t)}
{f(t)\over t-x} \kern 1pt dt,
\label{herm}
\end{equation}
where
\begin{equation}
\phi(z) = \prod_{k=0}^n(z-x_k)\biggl/\kern 1.5pt
\prod_{k=1}^n(z-p_k).
\label{phidef}
\end{equation}
\label{hermthm}
\end{theorem}
\medskip

\begin{figure}
\begin{center}
\vskip 10pt
\includegraphics[scale=.8]{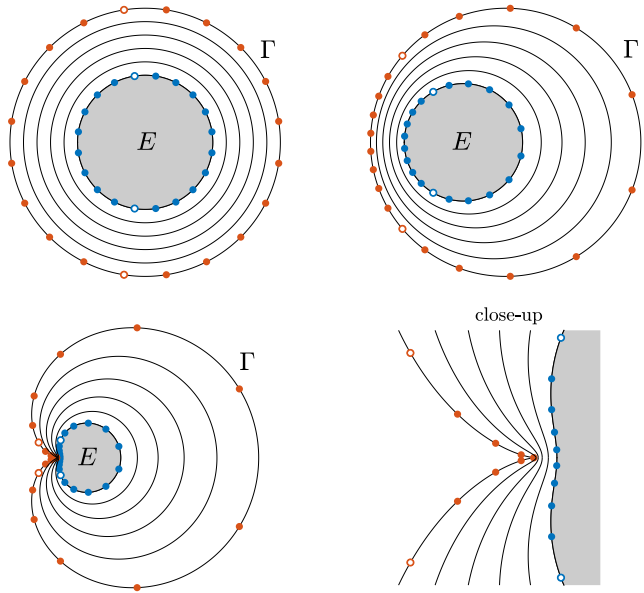}~~~~
\end{center}
\caption{\label{figpotent}Potential theory and
rational approximation.  In each image,
the shaded region is an approximation domain $E$ for a function
$f$ analytic in the region $\Omega$ bounded by
$\Gamma$.\ \ If we think of the poles of
an approximation $r\approx f$ as positive point charges
and the interpolation points as negative point charges,
then a minimal-energy equilibrium distribution of the charges
gives a favorable configuration for approximation.
This is a discrete problem of
potential theory that becomes continuous in
the limit $n\to \infty$, enabling one to take advantage of
invariance under conformal maps.
In these images $E$ and\/ $\Gamma$ are
disjoint and the convergence is exponential, but the third
domain and its close-up illustrate the clustering effect,
which will become more pronounced as the gap shrinks to zero.  The
pairs of interpolation points and poles marked by hollow dots delimit
one half of the total, highlighting how both sets of points
accumulate close to the singularity.}
\end{figure}

To see how this theorem is applied,
let\/ $\Omega$ be a simply connected domain bounded
by a closed curve~$\Gamma$,
as indicated in Fig.~\ref{figpotent} (see also Fig.~9 in
the next section),
and let $f$ be analytic in 
$\Omega$ and extend continuously to $\Gamma$.
Suppose $f$ is to be approximated on a
compact set $E\subset \Omega$, which in this section we
take to be disjoint from $\Gamma$.
Theorem~\ref{hermthm} implies that for any $x\in E$,
\begin{equation}
|f(x) - r(x)| \le C \tau,
\label{tauest}
\end{equation}
where $C$ is a constant independent of $n$ and $\tau$ is the ratio
\begin{equation}
\tau = {\max_{z\in E} |\phi(z)|\over
\min_{z\in \Gamma} |\phi(z)|}.
\label{rhodef}
\end{equation}
If $\phi$ is much smaller on $E$ than
on $\Gamma$, then $\tau$ and hence $f-r$ must be small.

Figure~\ref{figpotent} gives an idea of how this can happen.
In each image,
the red dots on~$\Gamma$ represent a good choice of poles $\{p_k\}$
and the blue dots on the boundary of
$E$ a corresponding good choice of interpolation points $\{x_k\}$.
Consider first the 
upper-left image, where $E$ and $\Gamma$ define a circular annulus. 
The equispaced configurations of $\{p_k\}$ and $\{x_k\}$
ensure that $\tau$ will decrease exponentially as $n\to\infty$.
To see this, in view of (\ref{phidef}), we define
\begin{equation}
u(z) = n^{-1}\sum_{k=0}^n \log|z-x_k| -
n^{-1}\sum_{k=1}^n\log |z-p_k|.
\end{equation}
This is the potential function generated by $n+1$ 
negative point charges of strength
$n^{-1}$ at the interpolation points and $n$ positive point
charges of strength $-n^{-1}$
at the poles.
Then $\exp(n u(z)) = |\phi(z)|$,
and therefore
\begin{equation}
\tau = \exp(-n \kern 1pt [\kern 1pt  \min_{z\in \Gamma} u(z) -
\max_{z\in E} u(z)\kern 1pt]\kern 1pt ).
\end{equation}
For $\tau$ to be small, we want $u$ to be uniformly bigger on
$\Gamma$ than on $E$.  Finding the best such configuration is an
extremal problem that will be approximately solved if the points
are placed in an energy-minimizing equilibrium position.  In each
of the images of Fig.~\ref{figpotent}, the points are close to
such an equilibrium.  Each charge is attracted to the charges of
the other color, but repelled by charges of its own color.

Finding an optimal configuration (for the given 
choice of $\Gamma$) is complicated for
finite~$n$, but the problem becomes cleaner
in the limit $n\to \infty$, and this is where
the power of potential theory is fully revealed.
We now imagine 
continua of interpolation points and poles defined by
a signed measure $\mu$ supported on $E$, where it is nonpositive
with total mass $-1$,
and on $\Gamma$, where it is nonnegative with total
mass $1$.  It can be shown that
there is a unique measure of this kind that minimizes the energy
\begin{equation}
I(\kern .7pt\mu) = - \int\kern -4pt \int \log|z-t| \kern .7pt d\mu(z) \kern .7pt d\mu(t),
\label{energy}
\end{equation}
with associated potential function
\begin{equation}
u(z) = -\int \log|z-t|\kern .7pt  d\mu(t),
\label{potential}
\end{equation}
and $u$ takes constant values $\uE < 0$ on $E$ and $\uG = 0$ on $\Gamma$.
The minimum $\Imin = \inf_\mu I(\kern .5pt \mu)$ is known
to be positive, and
for minimax degree $n$ rational approximations $r_n^*$ one has
exponential convergence as $n\to \infty$ at a corresponding rate:
\begin{equation}
\limsup \| f- r_n^*\|^{1/n} \le \exp(-\Imin).
\label{optrate}
\end{equation}
(The actual rate is in fact twice as fast as this, $\exp(-2\kern .3pt \Imin)$,
for functions whose singularities in the complex plane are
just isolated algebraic branch points~\cite[p.~93]{levsaff},
\cite{stahlgeneral}.)

The reciprocal of $\Imin$ is known as the {\em condenser capacity}
for the $(E,\Gamma)$ pair, a term that reflects an electrostatic
interpretation of the approximation problem.  In electronics,
capacitance is the ratio of charge to voltage difference. A
capacitor has high capacitance if its positive and negative plates
are close to one another, so that the attraction of charges of
opposite sign enables a great deal of charge to be accumulated
on them without the need for much of a voltage difference.
For fast-converging rational approximation, on the other hand,
we want $\Gamma$ and $E$ to be far apart, corresponding to a {\em
small\/} amount of charge relative to the voltage difference,
hence small numbers of poles and interpolation points needed to
achieve a given ratio $\tau$.

We can now see how the second and third images
of Fig.~\ref{figpotent} were drawn.  They were obtained by
conformal transplantation, exploiting the invariance of problems
of potential theory under conformal maps.  The eccentric domain
of the second image comes from a M\"obius transformation, and the
pinched domain of the third image comes from a further squaring.
The blue and red points obtained as conformal images of equispaced
points in the symmetric annulus are known as Fej\'er--Walsh
points~\cite{starke93}.

One might wonder, for arguments of this kind, is it necessary
to place the poles of\/ $r$ on the boundary of the region of
analyticity of $f\kern .7pt $?  In fact, $\Gamma$ does not
have to lie as far out as that boundary, nor do the poles have
to be on $\Gamma$, for as stated in Theorem~\ref{hermthm}, the
integral representation (\ref{herm}) is valid for any placement
of the poles.  Asymptotically as $n\to\infty$, however, it is
known that the convergence rate cannot be improved by placing
poles beyond the region of analyticity of $f$~\cite{levsaff}.  A
special choice is to put all the poles at $x=\infty$, in which
case rational approximation reduces to polynomial approximation,
still with exponential convergence though at a lower rate than
in~(\ref{optrate}).

\section{\label{sec5}Explanation of tapered exponential clustering}
Now we examine how the analysis of the last section
must change for approximations with singularities.
There is a considerable specialist literature here by
authors including Aptekarev, Saff, Stahl, Suetin, and
Totik~\cite{bt,levsaff,rtw,safftotik,stahl,stahl12,suetin}, which
investigates certain best approximations in detail.  Our emphasis
is on the broad ideas applicable to near-best approximations too.

\begin{figure}
\begin{center}
\vskip .2in
\includegraphics{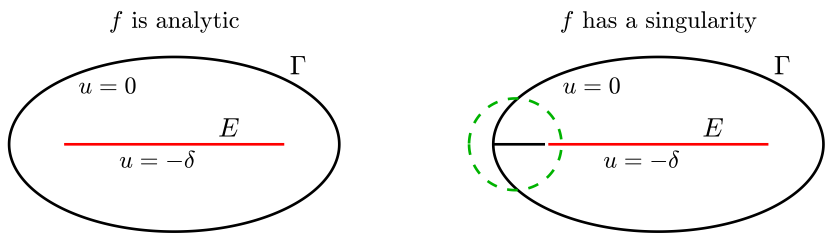}
\end{center}
\caption{\label{sketchfig}Two kinds of
problems of rational approximation of a function $f$ on a 
domain $E$.  On the left (section~\ref{sec4}), $f$ is analytic on $E$ and poles
can be placed on a contour\/ $\Gamma$ enclosing $E$ in the region of
analyticity:\ convergence is exponential with accuracy on the order
of $\exp(-n\kern .5pt \delta\kern .5pt)$ for a constant
$\delta>0$.
On the right (section~\ref{sec5}), $f$ has a singularity at a point $\zc$
on the boundary of $E$, and
$\Gamma$ must touch $E$ at $\zc$:\ convergence is root-exponential with
accuracy of order $\exp(-n\kern .5pt \delta\kern .5pt)$ again,
but now with~$\delta$ diminishing at the rate $1/\sqrt n$ as $n\to\infty$.
In the circled region, the potential makes the transition from
$\uG= 0$ to $\uE = -\delta$.}
\end{figure}

From Fig.~\ref{figpotent} it is clear that potential theory
should give some insight when $f$ has a singularity on
the boundary of $E$.  The lower pair of images shows clustering
of poles where $\Gamma$ has a cusp close to the boundary of $E$,
and as the cusp is brought closer to $E$, the clustering will
grow more pronounced.  However, the argument we have
presented breaks down when $\Gamma$ actually touches $E$.
The situation is sketched
in Fig.~\ref{sketchfig}.  Physically, this would be a capacitor
of infinite capacitance, implying that an equipotential distribution
$u$ with a nonzero voltage difference would require an infinite
quantity of charge.  Mathematically, the estimate (\ref{tauest})
fails because $\tau$ cannot be smaller than $1$.

\begin{figure}
\begin{center}
\vskip 10pt
\includegraphics[scale=.9]{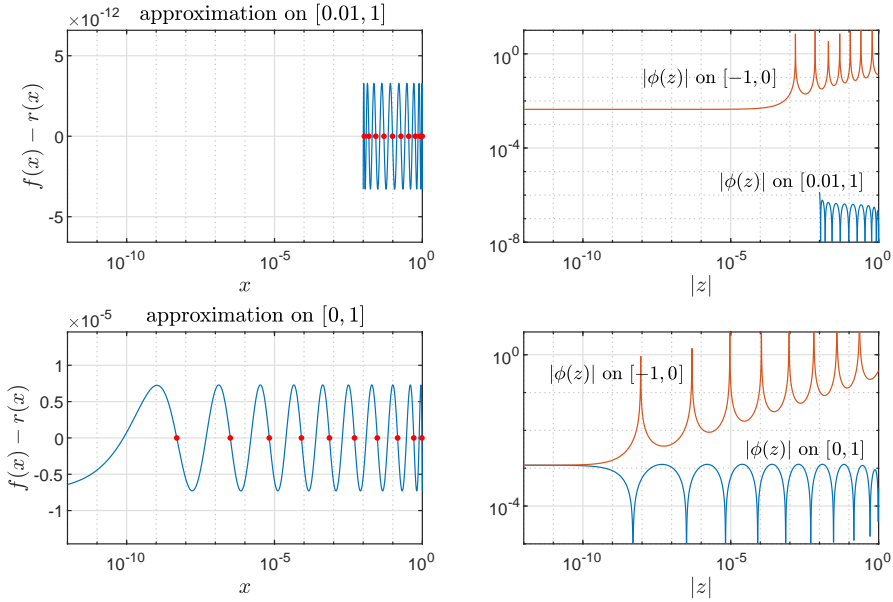}~~~~
\end{center}
\caption{\label{figmini}On the left, error curves in type $(9,10)$ 
minimax approximation of $\sqrt x$ on $\inzz$ and $\inz$.
On the right, plots of $\phi(z)$ as defined
by $(\ref{phidef})$ on these approximation intervals
and on $\inzm$.  The curves in the upper-right image show
a reasonable approximation to constant
values on $\inzm$ (upper curve) and on $\inzz$ (lower
curve), but in the lower-right image,
nothing like constant behavior of
$|\phi(z)|$ on $\inzm$ is evident.  We explain this
by noting that what matters to the accuracy of an approximation
is the integral (\ref{hermagain}) of $|\phi(x)/\phi(t)|$ with
respect to $t\in \Gamma$, not its maximum.  Taking
advantage of this property, poles and interpolation points
distribute themselves more sparsely near the singularity,
freeing more of them to contribute to the approximation further away---the
phenomenon of tapered exponential clustering.}
\end{figure}

To see what happens in such cases, we can examine
the function $\phi$ computed numerically for an example problem.
The left column of Fig.~\ref{figmini} shows error curves in type $(9,10)$
minimax approximation of $\sqrt x$ on $\inzz$ (above) and $\inz$
(below).  (Type $(m,n)$ means numerator degree at most $m$
and denominator degree at most $n\kern .3pt$; we choose these
parameters rather than $(n,n)$ to make the plots slightly cleaner.)
The curves each equioscillate between $m+n+2 = 21$ extrema,
and in the lower curve, on the semilogx scale, we see the
wavelength increasing as $x\to 0$.  As a minimax approximation
with free poles, this rational function has $m+n+1 = 20$ points
of interpolation rather than the standard number $n = 10$ for
an approximation with preassigned poles, so for the cleanest
display of the potential function $\phi$ in the right column we
have picked out just half of these, marked by the red dots.

The right column of Fig.~\ref{figmini} shows the function
$|\phi|$ plotted on the approximation interval
(the lower blue curve) and on the important portion $\inzm$
of the integration contour $\Gamma$ (the upper red curve).
(To be precise, for these plots the numerator of (\ref{phidef})
ranges over just the interpolation points $x_1,\dots, x_n$
marked by red dots.)  In the upper image, for $\inzz$, the curves
reveal a reasonable approximation to what the last section
has led us to expect from
potential theory.  The blue curve has approximately even magnitude,
and this is about five orders of magnitude below the red curve,
also of approximately even magnitude.  Thus the ratio $\tau$ of
(\ref{rhodef}) is far below $1$, and the estimate (\ref{tauest})
serves to bound the approximation error.  (The actual
error is about the square of this bound since we have omitted
half the interpolation points.)

The lower image, which is a centerpiece of
this paper, tells a strikingly different story.  Here again the
blue curve is flat, showing the even dependence on $x$ we expect
in a minimax approximation.  The red curve for $|\phi(z)|$ on
$\inzm$, however, is now tilted at an angle on these log-log axes,
showing a steady closing of the gap between the curves as\/ $t$
moves from $-1$ to $0$.
Clearly in this case $\inzm$ is not at all a curve of
constant $|\phi|$.

To understand the linearly closing gap in Fig.~\ref{figmini}, we note
that what fails in the analysis of the
last section for an approximation problem
with a singularity is not the Hermite integral,
\begin{equation}
f(x) - r(x) = {1\over 2\pi i} \int_\Gamma {\phi(x)\over \phi(t)}
{f(t)\over t-x} \kern 1pt dt,
\label{hermagain}
\end{equation}
but the
estimate (\ref{tauest}) we derived from it.  Implicitly (\ref{tauest})
came from bounding (\ref{hermagain}) by H\"older's inequality,
\begin{equation}
|f(x) - r(x) |  \le {1\over 2\pi}
\left\|{\phi(x)\over \phi(t)}\right\|_\infty 
\left\|{f(t)\over t-x}\right\|_\infty \|1\|_1^{},
\label{option1}
\end{equation}
where the $\infty$- and $1$- norms are defined over $t\in \Gamma$.
(The norm $\|1\|_1^{}$ is equal to the length of $\Gamma$.)
When $\Gamma$ and $E$ are disjoint, the first $\infty$-norm in (\ref{option1}) is
exponentially small as $n\to\infty$ and the second is bounded.
However, these properties fail as $\Gamma$ and $E$ touch.
We can rescue the argument
by noting that $|\phi(x)/\phi(t)|$ does not have to be
small for all $t$ so long as its integral is small.
More precisely, the quantity 
$f(t)/(t-x)$ of (\ref{option1}) may not be bounded as $t,x\to \zc$ but
$f(t)|\kern .5pt t-\zc|^{1-\alpha} /(t-x)$ will be
bounded if we assume
$f(t-\zc) = O(|\kern .5pt t-\zc|^\alpha)$ for some constant $\alpha$.
So what actually matters is that the integral of  
$|\kern .5pt t-\zc|^{\alpha - 1} |\phi(x)/\phi(t)|$
should be small, and we
accordingly replace (\ref{option1}) by the alternative H\"older estimate
\begin{equation}
|f(x) - r(x) |  \le {1\over 2\pi}
\left\|{\phi(x)\over \phi(t) }|\kern .5pt t-\zc|^{\alpha-1}\right\|_1
\left\|{f(t)\over t-x}|\kern .5pt t-\zc|^{1-\alpha}\right\|_\infty.
\label{option3}
\end{equation}

\begin{figure}
\begin{center}
\vskip 18pt
\includegraphics[scale=1.2]{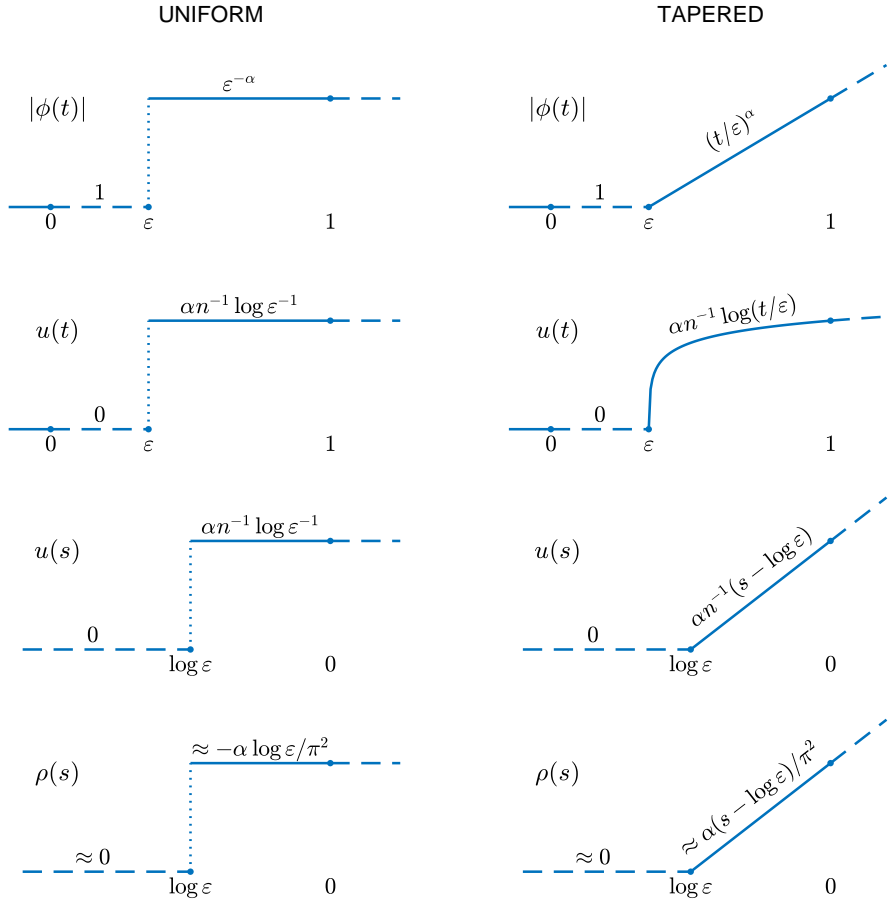}
\end{center}
\caption{\label{relupot}The potential theory of exponential clustering
(continuation of Fig.~{\rm \ref{relu}}).
The first two rows (right column) show the function $|\phi(t)|$ of $(\ref{phidef})$
and the associated potential $u(t) = n^{-1} \log |\phi(t)|$
for the model problem $(\ref{modelprob0})$--$(\ref{modelprob})$.
The third row shows the behavior along the real $s$-axis after
the change of variables to $s = \log t$; the domain is now
the infinite strip $0 < \Im s < \pi$, with $u=0$ for $\Im s = \pi$.
The final row shows the charge density $\rho(s) = n u_n(s)/\pi$, where
$u_n$ is the normal derivative of $u$ on the boundary of the strip.
The intervals that matter (emphasized by
solid rather than dashed lines)
are $\ve < |t| < 1$ in the $t$ variable
and $\log \ve < \Re s < 0$ in the $s$ variable.
Smaller values of $|t|$ and $s$ contribute
negligibly to the integral {\rm (\ref{hermagain})},
and larger values are far from the singularity.}
\end{figure}

For simplicity let us assume that the singularity
lies at $\zc = 0$ and the
part of the contour $\Gamma$ that matters is $\inz$, and
let the domain be scaled so that $|\phi(x)|\approx 1$ for $x\in E$.
We want $|\phi(t)|$ to be at most
$1$ for $t<0$ and at least $(t/\ve)^\alpha$ for $t> \ve$,
where $\ve$ is the distance of the closest pole to the singularity.
See the upper-right image of Fig.~\ref{relupot}.
Defining $u(t) = n^{-1} \log \phi(t)$ leads us
to the model problem sketched in the image below
this in the figure: find a harmonic function $u(t)$ in
the upper half $t$-plane such that
\begin{equation}
u(t) = \cases{
0, & $t \le \ve$, \cr\noalign{\vskip 3pt}
\alpha \kern .5pt n^{-1}\log(t/\ve) , & $t> \ve$.\cr
}
\label{modelprob0}
\end{equation}
We now make the change of variables $s = \log t$,
which transplants the Laplace problem
to the infinite strip $S = \{s\in\C:\, 0 <\Im s < \pi\}$, as
sketched in the $(3,2)$ position of the figure:
find a harmonic function $u(s)$ in $S$  satisfying 
\begin{equation}
u(s) = \cases{
0, & $\Im s = \pi$, \cr\noalign{\vskip 3pt}
0, & $\Im s = 0$ and $\Re s \le \log \ve$, \cr\noalign{\vskip 3pt}
\alpha \kern .5pt n^{-1}(s-\log\ve), & $\Im s = 0$ and $\Re s > \log \ve$.\cr
}
\label{modelprob}
\end{equation}
This change of variables is convenient mathematically, and it is
also important conceptually, since it is well known
that influences on harmonic functions decay
exponentially with distance along a strip.
Consequently, if $\ve$ is small, the solution to a Laplace problem
for $\log \ve \ll \Re s \ll 0$ will be essentially
(though not exactly) determined by the boundary conditions
in that region.
This just matches what we need for the model problem as posed in
the original~$t$ variable,
where behavior for $|t|$ of order
$\ve$ or less is unimportant because it contributes negligibly to the
integral (\ref{hermagain}) and behavior for $|t|$ of order
$1$ or more is unimportant because it is 
far from the singularity under investigation. 

So we address our attention to (\ref{modelprob}).
An exact solution can be obtained via the Poisson
integral formula for an infinite strip~\cite{widder},
\begin{equation}
u(x+iy) = {\alpha \kern .5pt n^{-1}\over 2\pi} \int_0^\infty {\xi \sin(y)\over
\,\cosh(\xi - (x-\log \ve)) - \cos(y)\,} \kern 1pt d\xi,
\label{solnexact}
\end{equation}
where we have set $s = x+iy$.  However, we do not need exactly this
since the region where our model applies is
$\log \ve \ll \Re s \ll 0$.  In this region, the bilinear harmonic function
\begin{equation}
u(x+iy) = \alpha\kern .5pt n^{-1} (1-{y\over \pi}) (x - \log\ve)
\label{solnapprox}
\end{equation}
satisfies the boundary conditions and is accordingly
a good approximation to the solution to~(\ref{modelprob}).  The corresponding
pole density distribution on the real $x$ axis is $(n/\pi)$ times
the normal derivative, 
\begin{equation}
\rho(s) = - {n\over \pi} {\partial \over \partial y} u(x+iy)
= {\alpha\over \pi^2} (x - \log\ve).
\end{equation}
This linear growth, sketched in the bottom-right image
of Fig.~\ref{relupot}, is just what we set out to explain in
Fig.~\ref{relu}.

Let us now look at the quantitative implications of this argument,
comparing uniform exponential clustering (left column of
Fig.~\ref{relupot}) with tapered exponential clustering
(right column).  
According to our model, the integral of the solid 
portions of the $\rho(s)$ curves in the bottom should be
equal to $n$, the total number of poles.
For uniform clustering the integral is $(-\alpha\kern .5pt\log \ve)^2/\pi^2$,
leading to the estimates
\begin{equation}
\hbox{Closest pole:~~} \ve \approx \exp(-\pi \sqrt{n/\alpha}\kern 1pt), \quad
\hbox{Accuracy:~~} \ve^\alpha \approx \exp(-\pi 
\sqrt{\alpha\kern .5pt n}\kern 1pt).
\label{est1}
\end{equation}
For tapered clustering the integral is ${1\over 2}(-\alpha\kern .5pt\log \ve)^2/\pi^2$, leading
to the estimates
\begin{equation}
\hbox{Closest pole:~~} \ve \approx \exp(-\pi \sqrt{2\kern .3pt n/\alpha}\kern 1pt), \quad
\hbox{Accuracy:~~} \ve^\alpha \approx \exp(-\pi 
\sqrt{2\kern .3pt \alpha\kern .5pt n}\kern 1pt).
\label{est2}
\end{equation}
Thus, as mentioned in the abstract, our model
leads to the prediction of a
factor of 2 speedup with tapered clustering.
It would be interesting to investigate whether, for certain problems,
exactly this ratio could be established theoretically
in the limit $n\to\infty$.

\begin{figure}
\vskip 1.4em
\begin{center}
\includegraphics[scale=.99]{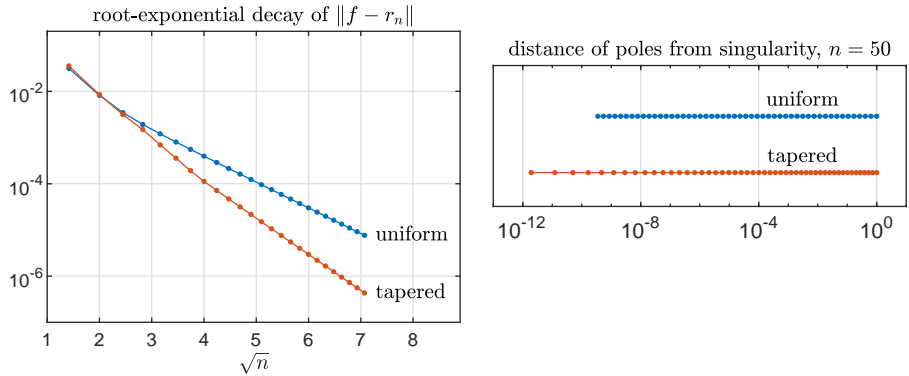}
\end{center}
\caption{\label{fig12}Linear-minimax approximation of $f(x) =\sqrt x$
on $\inz$ with preassigned exponentially
clustered poles in $\inzm$, $n = 2,4,\dots, 50$.
Tapering the distribution makes 
the convergence rate approximately double, as predicted by the model of 
Section~$\ref{sec5}$.}
\end{figure}

As an example of a
problem in which we may make such a comparison
numerically, consider Fig.~\ref{fig12}.
These data show $\infty$-norm errors for rational linear-minimax
approximations of even degrees $n$ from 2 to 50 with
preassigned exponentially clustered poles.  That is, the
approximations are optimal in the $\infty$-norm among rational 
functions in $\Rn$ with simple poles at the prescribed points; they
are characterized by error curves equioscillating
between $n+2$ extrema.
The upper curves
correspond to uniformly clustered poles
$p_k = -\exp(-\pi k/\sqrt n)$, $0\le k \le n-1$, and the lower curves
to tapered poles
$p_k = -\exp(\sqrt 2 \pi(\sqrt k-\sqrt n))$, $1\le k \le n$.
The asymptotic errors appear to be about
$\exp(-\sqrt{2.3 n})$ for uniform clustering and
$\exp(-\sqrt{4.7 n})$ for tapered clustering.
With $\alpha = 1/2$ for $f(x) = \sqrt x$, the corresponding
estimates (\ref{est1}) and (\ref{est2}) are
$\exp(-\sqrt{2.2 n})$ and $\exp(-\sqrt{4.4 n})$.

Analyses related to the argument we have presented were published
by Stahl for rational minimax approximation of $|x|$
on $[-1,1]$ and $x^\alpha$ on
$\inz$~\cite{stahl92,stahl93,stahl94,stahl}.
For $x^\alpha$ Stahl gives the result
\begin{equation}
\hbox{Accuracy:~~} \ve^\alpha \approx \exp(-\pi 
\sqrt{4\kern .3pt \alpha\kern .5pt n}\kern 1pt),
\label{stahlest2}
\end{equation}
which is not just an estimate but a theorem concerning
the limit $n\to\infty$ 
(assuming $\alpha$ is not an integer), with precise constants.
This is exactly what one would expect based on (\ref{est2}),
since, as mentioned earlier,
the effective value of $n$ is doubled
in the case of true minimax approximants~\cite{rakh}.

Stahl worked
essentially in the variable $t$ rather
than $s$, so his boundary conditions involved logarithms, as
in the second image of the right column of Fig.~\ref{relupot}.
Whenever one has a Laplace problem with Dirichlet boundary
data, one can interpret it as the problem of finding an
equipotential distribution in the presence of an external
field defined by that boundary data, and this interpretation
has been carried far in approximation theory~\cite{safftotik}.
From this point of view one can say that tapered exponential
clustering results from poles and zeros being slightly pushed
away from a singular point by a logarithmic potential field.

\section{\label{sec6}Exponential and double exponential quadrature}
In this final section we turn to another problem where exponential
clustering appears.
Let $f$ be a continuous function on $[-1,1]$.
We wish to approximate the integral of $f$ by a linear combination
\begin{equation}
I_n = \sum_{k=1}^n w_k f(x_k),
\label{quadsum}
\end{equation}
where $\{x_k\}$ are distinct nodes in $[-1,1]$ and
$\{w_k\}$ are corresponding weights, in such a way that
the accuracy is good even if $f$ has
branch point singularities at the endpoints.
To this end, we introduce a change of variables $g(s)$ from the real
line to $[-1,1]$, so that the integral becomes
\begin{equation}
I = \int_{-1}^1 f(x) \kern .7pt dx = \int_{-\infty}^\infty f(\kern .5pt g(s))
\kern .7pt g'(s) \kern .7pt ds,
\label{cov}
\end{equation}
and we apply the equispaced trapezoidal rule.  This involves
an infinity of sample points in principle, but if $g'(s)$ decays
rapidly, we may truncate
these to an $n$-point rule like (\ref{traprule}):
\begin{equation}
I_n = \kern 2pt h \kern -6pt \sum_{k =-(n-1)/2}^{(n-1)/2}
f(\kern .3pt g(kh)) \kern .7pt g'(kh).
\label{DESE}
\end{equation}
Quadrature formulas of this kind were introduced 
around 1970 by Mori, Takahasi,
and other Japanese researchers and also in
the analysis of sinc methods by Stenger.
See~\cite{haber,IMT,stengersurvey,stengerbook,tm73,trap},
as well as~\cite{mori} for the history as told by Mori himself.
The standard ``exponential'' choice of $g$ is
\begin{equation}
g(s) = \tanh(s), \quad g'(s) = \sech^2(s),
\label{tanh}
\end{equation}
with which (\ref{DESE}) becomes the {\em tanh formula}.
As in section 2, we estimate the truncation error as of order
$\exp(-nh)$ and the discretization error of order $\exp(-\pi^2/h)$.
(The latter could be worse if $f$ has additional singularities
near $(-1,1)$.)  This gives a balance $h \approx \pi/\sqrt n$,
with convergence rate of order $\exp(-\pi\sqrt n\kern .7pt )$.
An estimate of this form is valid for any H\"older continuous
branch point singularity; see~\cite[Thm.~3.4]{stengersurvey},
\cite[Thm.~2.1]{tanaka}, and~\cite[Thm.~14.1]{trap}.

Root-exponential convergence!  This is much better than any
algebraic order, but for practical applications on one-dimensional
domains, methods of this kind often seem very wasteful, with almost
all the points
being used up in resolving the singularity
(100$\%$ of them, in the limit $n\to \infty$)~\cite{sincfun}.
A year or two after the first exponential formulas
appeared, it was realized that one can do better with
``double exponential'' formulas.
We focus on the {\em tanh-sinh formula\/} proposed by Takahasi
and Mori in~\cite{tm74} and subsequently used and analyzed
by many others including Okayama, Sugihara, and Tanaka as
well as Bailey and Borwein~\cite{bb,ms01,oms,sugihara,tanaka}.
Here (\ref{tanh}) is replaced by
\begin{equation}
g(s) = \tanh(\pih \sinh(s)), \quad g'(s) =
      \pih\cosh(s)\kern 1.3pt \sech^2(\pih\sinh(s)).
\label{tanhsinh}
\end{equation}
Under suitable assumptions we can now estimate the truncation 
and discretization errors
as of orders $\exp(-(\pi/2) \exp(nh/2))$ and
$\exp(-\pi^2/h)$.  The first of these estimates is
the big improvement, for this quantity can be almost-exponentially
small with a much smaller value of $h$ than before, of order
$\log(n)/ n$ rather than $1/\sqrt n$.
By {\em almost-exponential,} we mean of order\break
$\exp(-C n/\log n)$ for some \hbox{$C>0$}.
With this reduced value of $h$, the second estimate
becomes almost-exponentially small too.

\begin{figure}
\begin{center}
\vspace{12pt}
\includegraphics[scale=.96]{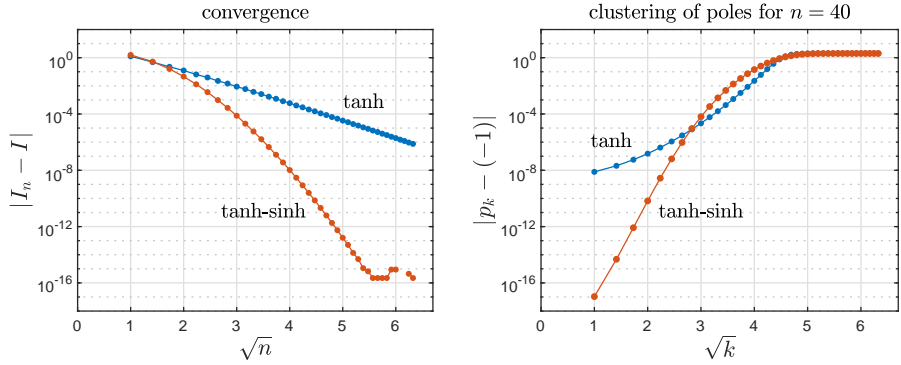}
\vspace{-7pt}
\end{center}
\caption{\label{fig13}
On the left, root-exponential convergence of the tanh quadrature
formula applied to integration of $\sqrt{1+x}$ (note the $\sqrt n$
axis as usual); the tanh-sinh formula converges much faster down
to machine precision.
On the right, the distances of nodes from poles (with a $\sqrt k$
axis) show uniform exponential clustering for the
tanh formula with $n=40$ and tapered exponential clustering for tanh-sinh.}

\end{figure}

Figure~\ref{fig13} shows data for the tanh and tanh-sinh formulas.
(We used the empirical choices $h = \pi/\sqrt n$ and $h = 1.2\log(2\pi n)/n$,
respectively.)
The left image plots $|I_n- I|$  against $\sqrt n$
for $n$ from $1$ to $40$ for the integration
of $f(x) = \sqrt{1+x}$.  The tanh curve appears
straight, confirming the root-exponential convergence, and
the tanh-sinh curve bends downward, confirming that
its rate is faster.  The unexpected image is on the right, a plot of
distances of the nodes from the endpoint $x=-1$.
For tanh quadrature, these distances are uniformly exponentially spaced,
appearing as a parabola on these axes.
The curve for tanh-sinh quadrature, however, is
almost perfectly straight.  It would seem that tanh-sinh quadrature
exploits tapered exponential clustering!
It surprised us when we first saw curves like this.
Why is there a resemblance between the tanh and tanh-sinh
quadrature formulas and the phenomena of rational approximation
discussed in the earlier sections of this article?

Some steps toward an answer come from a beautiful connection
introduced by Gauss and exploited by Takahasi and 
Mori~\cite{gauss,tm71,gaussCC}: every quadrature formula
can be associated with a rational approximation.
Suppose first that $f$ can be analytically
continued to a neighborhood $\Omega$ of $[-1,1]$ bounded by a
contour $\Gamma$.
Then the integral can be written
\begin{equation}
I = \int_{-1}^1 f(x) \kern .7pt dx = {1\over 2\pi i} \int_\Gamma
f(t) \kern .7pt \varphi(t)\kern .7pt  dt,
\label{Iint}
\end{equation}
where the {\em characteristic function} $\varphi$ is defined by
\begin{equation}
\varphi(t) = \int_{-1}^1 {dx\over t-x} = \log{t+1\over t-1}.
\end{equation}
On the other hand the quadrature sum (\ref{quadsum}) can be written
\begin{equation}
I_n = {1\over 2\pi i} \int_\Gamma f(t) \kern .7pt r(t)\kern .7pt  dt, 
\label{Inint}
\end{equation}
where $r$ is the degree $n$ rational function defined by
\begin{equation}
r(t) = \sum_{k=1}^n {w_k\over t-x_k}.
\label{ratdef}
\end{equation}
Subtracting (\ref{Inint}) from (\ref{Iint}) gives what we
call the {\em Gauss--Takahasi--Mori (GTM) contour integral,}
\begin{equation}
I - I_n = {1\over 2\pi i} \int_\Gamma^{} f(t) \kern .7pt
(\varphi(t)-r(t)) \kern .7pt dt
\label{rel1}
\end{equation}
and the corresponding error bound
\begin{equation}
| I - I_n | \le {1\over 2\pi}
\|f\|_\infty^{}\kern .7pt \|\varphi - r\|_\infty^{}\kern .7pt \|1\|_1^{},
\label{rel2}
\end{equation}
which we have written in the style of (\ref{option1}),
with the norms defined over $\Gamma$.

Equations (\ref{rel1}) and (\ref{rel2}) relate accuracy of a
quadrature formula to an approximation problem: if the nodes
and weights are such that $\varphi - r$ is small on the boundary
$\Gamma$ of a region where $f$ is analytic, then $|I-I_n|$ must
be small.  This reasoning was applied by Takahasi and Mori to
a range of quadrature formulas~\cite{tm71}.  Now $\varphi$ is an
analytic function in the extended complex plane minus the segment
$[-1,1]$.  It follows that so long as $\Gamma$ is disjoint from
$[-1,1]$, there exist rational approximations to $\varphi$ that
converge exponentially on $\Gamma$ as $n\to\infty$.  In particular,
this holds for the rational functions associated with Gauss and
Clenshaw--Curtis quadrature~\cite{gaussCC}, where it is convenient
to take $\Gamma$ in the form of an ellipse about $[-1,1]$ with foci
$\pm 1$.  It follows that both these quadrature formulas
converge exponentially as $n\to\infty$ for analytic integrands
(cf.\ \cite[Thm.~19.3]{atap}).

But what if $f$ has endpoint singularities?  Now $\Gamma$ must
touch $[-1,1]$ at the endpoints, and (\ref{rel2}) fails just
as (\ref{option1}) did in such a case.  In fact, this failure
is more severe, since $\|\varphi - r\|_\infty = \infty$
for any $r$ because of the logarithmic singularities of
$\varphi$.
The last section, however, suggests a solution.
Instead of (\ref{rel2}), we can derive from (\ref{rel1}) the bound 
\begin{equation}
| I - I_n | \le {1\over 2\pi} \kern .7pt\|f\|_\infty{}
\|\varphi - r\|_1^{} \kern 1pt .
\label{rel3}
\end{equation}
The switch from the $\infty$- to the $1$-norm changes the
problem of rational approximation of $\varphi$ profoundly.
Since the dominant effects just concern approximation of a logarithmic
singularity near the singular point, the
essential question becomes, how fast can $\log t$ be approximated
by rational functions over both sides of the interval $\inzm$ in
the $1$-norm? 

\begin{figure}
\begin{center}
\vspace{12pt}
\includegraphics[scale=1]{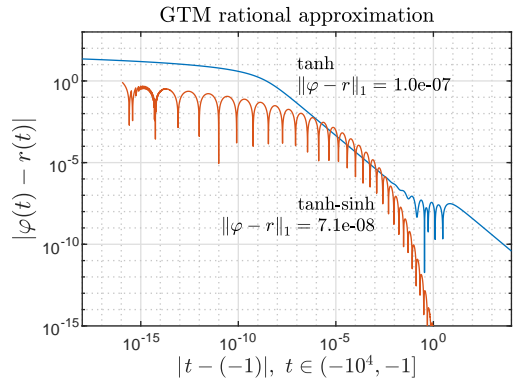}
\vspace{-7pt}
\end{center}
\caption{\label{fig14} Error $|\varphi(t)-r(t)|$ as a function of
distance to the left from $t=-1$ for the tanh and tanh-sinh
approximations with $n=40$ with 
$h = \pi/\sqrt n$ and $h = 1.2\log(2\pi n)/n$, respectively.
By symmetry, the same behavior would appear to the right from $t=1$.
Compare Fig.~\ref{figmini}, where the ratio of the blue and
red values in the lower-right image is closely analogous 
to the blue curve here.
The $1$-norms of the approximation errors over $[-2,-1]$ are indicated.
The slight irregularities at the left are the result of rounding error.}
\end{figure}
 
As we did with Fig.~\ref{figmini}, let us get some insight
by looking at the details of the approximation problem.
The rational function (\ref{ratdef}) for the tanh rule is
\begin{equation}
r(t) = \kern 2pt h \kern -6pt
\sum_{k =-(n-1)/2}^{(n-1)/2} {\sech^2(kh) \over t - \tanh(kh)},
\label{tanhrule}
\end{equation}
and for the tanh-sinh rule, it is
\begin{equation}
r(t) = \kern 2pt h \kern -6pt
\sum_{k =-(n-1)/2}^{(n-1)/2} {\pih\kern -.5pt  \cosh(kh)
\kern 1.3pt\sech^2(\pih\sinh(kh))
\over t - \tanh(\pih\sinh(kh))}.
\label{tanhsinhrule}
\end{equation}
Figure~\ref{fig14} plots $|\kern .5pt \varphi(t)-r(t)|$ for these
two approximations.

For tanh quadrature, we know that $|\kern .5pt
\varphi(t)-r(t)|$ must diverge to $\infty$ as $t\to -1$ because of
the log singularity of $\varphi$ at $t=-1$.  Yet the singularity is
so weak that the divergence only shows up as a gentle upward drift
in the blue curve at the left.  Over the main part of the plot, $|\kern
.5pt\varphi(t)-r(t)|$ decreases steadily down to around
$10^{-7}$.  The $1$-norm, measured here
over $t\in [-2,-1]$, is consequently very small, confirming via
(\ref{rel3}) the high accuracy of this quadrature rule.
As $n\to\infty$, this $1$-norm decays root-exponentially.

For tanh-sinh quadrature, again no approximation of $\varphi$
is possible in the $\infty$-norm.
In the $1$-norm, however, one might expect that the convergence
will now be almost-exponential.
Indeed, $|\kern .5pt \varphi - r|$ decays almost-exponentially as $n\to\infty$
over any domain bounded away from the singularity.  
But the 1-norm decay over the whole interval is in fact
just root-exponential, as is suggested by the number listed
being barely smaller than before.
The following reasoning suggests why this must be.
Consider approximation of $f(x) = \log x$ on $[\kern .5pt 0,
1]$.  Suppose rational approximations existed with faster than
root-exponential convergence in the $1$-norm.  Then by integrating,
we would get rational approximations to $g(x) = x\log x - x$ with
faster than root-exponential convergence in the $\infty$-norm,
which would contradict the evidence of Fig.~\ref{fig2}.

If $\|\kern .7pt \varphi - r\kern .7pt \|_1$ decreases only
root-exponentially as $n\to\infty$, how does the quadrature
formula converge almost-exponentially?  It appears that this
depends on additional properties
that go beyond rational approximation, involving
analytic continuation of the integrand onto an infinitely-sheeted
Riemann surface in exponentially small neighborhoods
of the endpoints~\cite{sugihara,tanaka}.

There remains the phenomenon of tapered exponential clustering, so
vividly evident in Fig.~\ref{fig13}.  We do not yet have
an explanation for this, nor a view of whether an
approximate $\sqrt k$ dependence
is genuine or just an artifact.  
This is a topic for ongoing research, where it would
be good to investigate also the 
distributions of exponentially clustered nodes, also
apparently tapered, that arise with the
``universal quadrature'' formulas of Bremer, et al.~\cite{brs,serkh}.

\section{Discussion}
\label{sec-discussion}
Exponential clustering of poles at singularities has been part
of the landscape of rational approximation for half a century,
but we believe this is the first study to focus on this effect.
Our motivation is that this clustering is what makes rational
approximations so powerful, and understanding it enables one
to improve existing numerical algorithms and develop new ones.
We find these phenomena fascinating,
especially the tapered clustering effect, and discovering that
tapering also appears in double exponential quadrature
was a bonus.  The elucidation of these matters with
the help of a sometimes seemingly endless program of numerical
experiments will forever be associated in our minds with the
Covid-19 shutdowns of 2020.

Here are some details of our computations.
Figures~\ref{fig1}, \ref{fig2}, \ref{fig6} and~\ref{figmini}
made use of the Chebfun {\tt minimax} command~\cite{minimax}, 
principally due to Silviu Filip, and Filip also kindly provided
us with a modified code for the weighted
minimax approximations of Figs.~\ref{fig2} and~\ref{fig6}.\ \ For
successful results in some of these
problems, we applied a M\"obius
transformation of $\inz$ to itself to weaken the singularity
while preserving the space $\Rn$.
For the approximations of Figs.~\ref{fig2}
and~\ref{fig6} on a complex disk, the AAA-Lawson algorithm was used as
implemented in Chebfun~\cite{chebfun,lawson}, again with a
M\"obius transformation.  Figure~\ref{fig3} was produced with
the {\tt confmap} code available at~\cite{laplace}, which in
turn calls {\tt aaa} from Chebfun~\cite{aaa} and {\tt laplace}
from~\cite{laplace}.
The {\tt aaa} code was also used directly
in Figs.~\ref{fig1} and~\ref{fig6},
and {\tt laplace} in Fig.~\ref{fig5}.
The Stokes and Helmholtz results of Fig.~\ref{fig5} were
produced by experimental codes that are not yet
publicly available developed with Abi Gopal and
Pablo Brubeck, respectively.
In Fig.~\ref{fig12}, a least-squares problem was extended by
a Lawson iteration (iteratively reweighted least-squares) to
compute minimax approximations with preassigned poles.
All the remaining results are based on straightforward
computations in \MATLAB\ and Chebfun.

\begin{acknowledgements}
We have benefited from helpful
advice from Bernd Beckermann, Pablo Brubeck,
Silviu Filip, Abi Gopal, Stefan G\"uttel, Arno Kuijlaars,
Andrei Mart\'inez-Finkelshtein,
Ed Saff, Kirill Serkh, Alex Townsend, and Heather Wilber.
\end{acknowledgements}

\end{document}